\documentclass[11pt]{article}
\usepackage[top=1in,bottom=1in,left=1in,right=1in]{geometry}

\usepackage{amsmath,amssymb,amsthm,mathrsfs}
\usepackage{algpseudocode}
\usepackage{multirow}
\usepackage{graphicx}
\usepackage{enumitem}
\usepackage{subcaption}
\usepackage{comment}
\usepackage[hyphens]{url}

\RequirePackage{natbib}
\usepackage{natbib}
\newcommand{\longPaperTitle}{School Bus Routing Problem with Open Offer Policy: incentive pricing strategy for students that opt out using school bus}

\usepackage{hyperref}
\hypersetup{
	unicode=false,
	pdftoolbar=true,
	pdfmenubar=true,
	pdffitwindow=false,
	pdfstartview={FitH},
	pdftitle={\longPaperTitle},
	pdfauthor={Caceres, Duran, Lespay, Contreras and Batta},
	pdfsubject={Operations Research},
	pdfnewwindow=true,
	bookmarksopen=true,
	colorlinks=true,
	linkcolor=red,
	citecolor=blue,
	filecolor=magenta,
	urlcolor=cyan
}

\usepackage{authblk}

\title{\longPaperTitle}

\author[1]{Hernan Caceres\thanks{hcaceres@ucn.cl (corresponding author)}}
\author[1]{Macarena Duran\thanks{mdp006@alumnos.ucn.cl}}
\author[1]{Hernán Lespay\thanks{hernan.lespay@ucn.cl}}
\author[1]{Juan Pablo Contreras\thanks{juan.contreras01@ucn.cl}}
\author[2]{Rajan Batta\thanks{batta@buffalo.edu}}
\affil[1]{Departamento de Ingeniería Industrial, Universidad Católica del Norte, Antofagasta, Chile}
\affil[2]{Department of Industrial and System Engineering, University at Buffalo, New York, USA}
\providecommand{\keywords}[1]{\textbf{Keywords:} #1}

\usepackage[colorinlistoftodos]{todonotes}

\usepackage{setspace}
\onehalfspace

\newcounter{problemno}
\DeclareRobustCommand{\newProblem}{%
	\mathscr{P}_\theproblemno}
\newcommand{\problemref}[1]{\hyperref[#1]{$\mathscr{P}_{\ref{#1}}$}}

\newtheorem{proposition}{Proposition}

\theoremstyle{definition}

\begin{document}


\maketitle

\begin{abstract}
This paper introduces the School Bus Routing Problem with Open Offer Policy (SBRP-OOP) that seeks to improve capacity usage and minimize the bus fleet by openly offering a monetary incentive to students willing to opt out of using a bus. We propose a mathematical formulation to determine a pricing strategy that balances the trade-off between incentive payments for students who choose not to use the bus with the expected savings obtained from operating fewer buses. To evaluate the effectiveness of the approach, we conducted simulations using both synthetic and real instances from a real operational context in the Williamsville Central School District (WCSD) of New York.
 
\end{abstract}


\section{Introduction}
School Bus Routing Problem (SBRP) is part of the family of Vehicle Routing Problem (VRP), in which a set of school bus routes has to be determined in order to pick up a given set of students from a set of potential bus stops and transport them to the school.  SBRP is cost-intensive for school administrators \citep{Ellegood2020}; therefore, improving transportation efficiency could help better manage and redistribute school resources, which could be used to improve educational programs, build new infrastructure, increase teachers’ wages, among others.

Typically, school districts assign all students to a bus route by default, regardless of whether they actually use the bus. Students are not incentivized to actively opt out when they do not use the transportation system. However, if students were allowed to opt out at the beginning of the school year, motivated by a monetary incentive, this could allow the district to design a more efficient and smaller set of routes exclusively for those students who choose to use the bus.


In this paper, we introduce the School Bus Routing Problem with Open Offer Policy (SBRP-OOP). The SBRP-OOP seeks to determine a pricing strategy to balance incentive payments to the students willing to opt out of using a bus and the savings obtained from using fewer buses to pick up students from their assigned stops. Students who accept the incentive are not considered in the routing generation process. Then, we consider an SBRP with a single school whose objective is to simultaneously: i) ﬁnd the set of stops to visit from a set of potential stops, ii) assign each student to a stop, considering a maximum total distance that a student has to walk and a maximum number of students per stop, and iii) generate routes to visit the selected stops; such that the total distance traveled by buses is minimized.

\subsection{Problem motivation}

\begin{figure}
	\begin{center}
		\includegraphics[width=0.8\columnwidth]{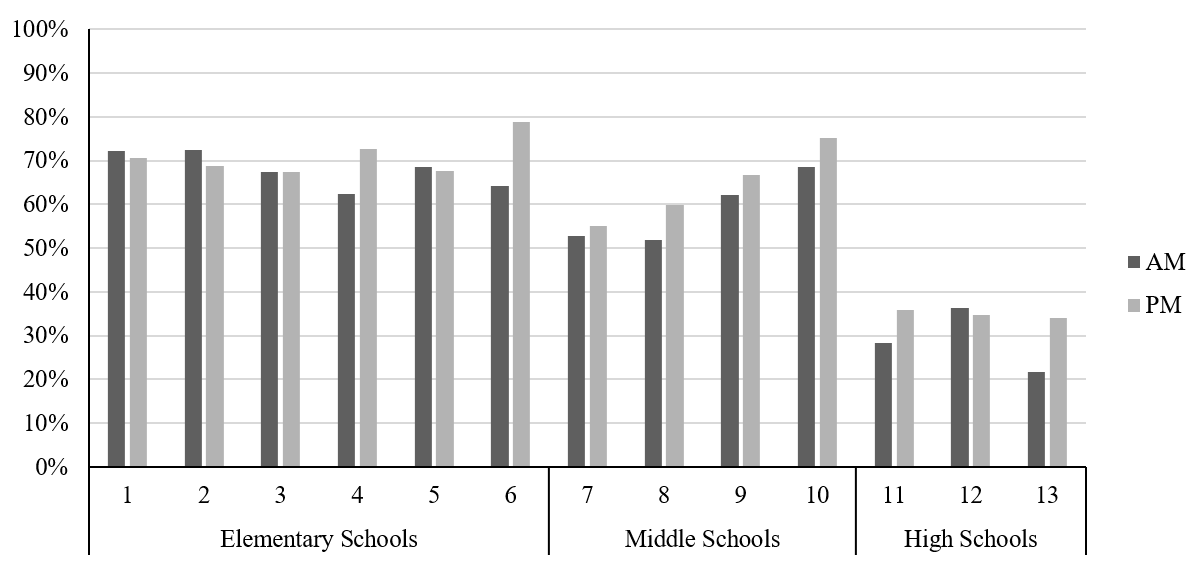}
		\caption{Average ridership per school and level.} \label{tax:ridership}
	\end{center}
\end{figure}

The SBRP-OOP is motivated by the operating context of the Williamsville Central School District (WCSD), which belongs to the New York State Education Department. School bus transportation in New York State is free and must be provided to all students by law. WCSD policy states that all students must be assigned to a stop and that all stops must be visited despite the uncertainty that students do not show up. This operational framework results in a dependency on a sizable bus fleet, with a substantial portion being underutilized. Consequently, the district faces considerable expenses related to bus maintenance and driver salaries, which leads to excessive costs. Figure \ref{tax:ridership} shows the ridership per school, separating the AM and PM cases. It can be observed that no more than 80\% of the students use the school transportation provided.

In addition to under-use, the WCSD is grappling with a severe shortage of bus drivers, a problem that the COVID-19 pandemic has exacerbated, as reported by the local media \citep{MikeDesmond2022}. Prospective candidates, typically retirees, harbor skepticism about working close to children and often prefer work-from-home or package delivery jobs \citep{lando2023}. As a result, the district faces significant challenges in both fleet reduction to achieve cost savings and ensuring an adequate number of personnel to meet the demands of its bus transportation system.



The WCSD uses an overbooking policy to improve bus utilization. However, the problem of idle capacity remains, which motivates exploration of new and potentially more risky solutions. The open offer policy aims to offer an incentive openly, extended as a check or a tax return, to any student willing to opt out of using a bus. The incentive amount is determined before the students decide whether to accept the deal and is the same for every student. Thus, the incentive value must be determined in such a way that it should reduce the risk of having to pay more for incentives than what would be saved by using fewer buses. 

\begin{figure}
	\begin{center}
 \includegraphics[width=1.0\columnwidth]{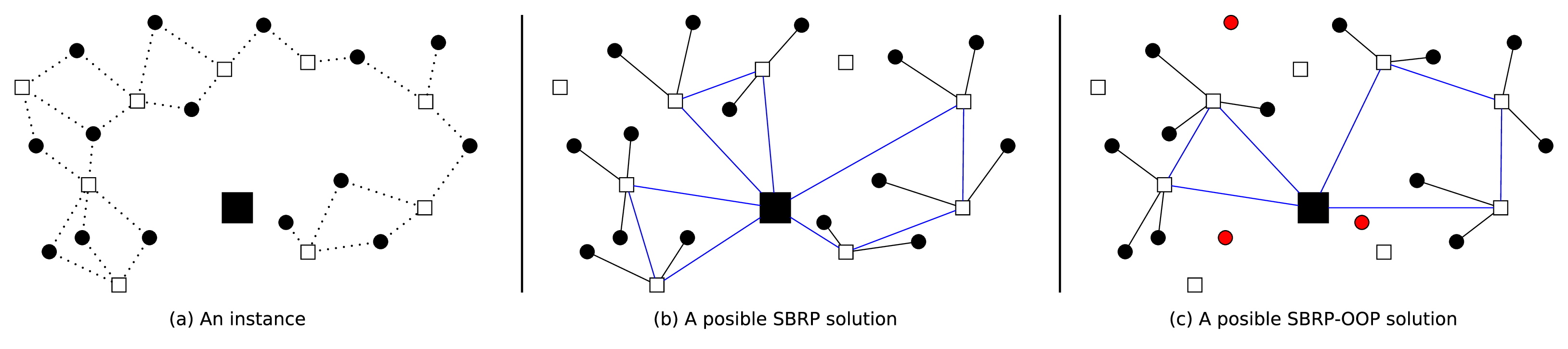}
		\caption{Motivation example for the SBRP-OOP} \label{tax:Motivation}
	\end{center}
\end{figure}

Figure \ref{tax:Motivation} illustrates the fleet reduction. Figure \ref{tax:Motivation}  (a) describes an example of a problem with a single school (black square), a set of potential stops (white square), and the set of students (black dots) connected to the stops that they are able to reach. Figure \ref{tax:Motivation} (b) shows a solution for the traditional SBRP with selection of bus stops (SBRP-BSS) when the bus capacity is six students. This SBRP solution uses three routes (represented by the blue lines), visiting seven stops. Finally, Figure \ref{tax:Motivation} (c) describes a solution of SBRP-OOP for the same bus capacity. In this case, three students accepted the incentive to opt out of using the buses (red dots). Therefore, these three students are not considered in the route planning process, and the fleet size is reduced to two buses that visit five stops. 

\subsection{Our contribution}

The contributions of this paper are fourfold.

First, we introduce the SBRP with an open offer policy and define performance indicators to measure the effectiveness of incentives in terms of expected savings and probability of failure. Second, we propsoe a new stop selection and routing generation model that uses a lexicographic objective function to first determine the minimum number of buses and then minimize the overall travel time of the buses. Our model considers classical routing and allocation constraints, as well as overbooking constraints based on the individual level of ridership of the students. Third, we propose a simulation-based solution approach to estimate expected savings and the probability of failure associated with an incentive value. To achieve this, we propose methods for estimating individual ridership and the opt-out probability of students based on average ridership, the amount of incentive, and the distance to the school. Then, we simulate the opt-out process and determine the minimum fleet required for students who decline the economic incentive. Finally, we conduct numerical experiments using both synthetic and real instances of the WCSD. We discover that the open offer policy can be applicable for certain schools saving up to three buses per school.

\subsection{Roadmap}

The remainder of this paper is organized as follows. Section~\ref{literature} gives a review of the literature of SBRP with bus stop selection and pricing demand management for routing problems. Section~\ref{formulation} presents a mathematical program for the SBRP with overbooking. Section~\ref{solution} describes our simulation-based methodology to tackle SBRP with open-offer policy. Section~\ref{experiments} reports the results of computational experiments with synthetic data, and Section~\ref{sec:casestudy} presents our results for real data from the WCSD. Finally, Section~\ref{discussion} presents the conclusions.

\section{Literature review} \label{literature}

The School Bus Routing Problem (SBRP) was introduced in \cite{newton1969} and has been extensively studied. SBRP consists of four interrelated subproblems \citep{Park2010}: the tactical bus stop selection subproblem, the operational bus route generation, bus route scheduling, and school bell time adjustment subproblems. For a review of the different classifications, objectives and contemporary trends in SBRP, we recommend the reader to \cite{Park2010, Ellegood2020}. The literature shows that the latest studies emphasize real-world issues such as heterogeneous fleets, mixed loading, multiple schools, and uncertainty of ridership \citep{Ellegood2020}. However, to our knowledge, demand management decisions focused on pricing strategies in which students are offered a monetary incentive to opt out of using the school transportation system have not received attention in the SBRP literature.

We present the relevant literature associated with the main characteristics of our studied problem: i) SBRP with bus stop selection and ii) pricing demand management for routing problems.

\subsection{SBRP with bus stop selection}
Because of the complexity of solving the bus stop selection and bus route generation subproblems simultaneously, these problems have been solved sequentially. Heuristic solution approaches can be classified into the location-allocation-routing (LAR) strategy and the allocation-routing-location (ARL) strategy (\cite{Park2010}). 

The LAR strategy determines a set of bus stops for a school and allocates students to these stops. Then, routes are generated for the selected stops. The main drawback of this strategy is that it tends to create excessive routes because restrictions on vehicle capacity are ignored in the location-allocation phase (\cite{Bowerman1995}). The main articles that use a heuristic approach based on the LAR strategy are \cite{bodin1979, dulac1980, li2021, farzadnia2021, calvete2021}. In \cite{bodin1979} the bus route scheduling subproblem with mixed loads is additionally considered. In \cite{dulac1980} multiple schools are considered, but mixed loads are not allowed; therefore, the problem is solved for each school independently. In \cite{Parvasi2017} the student's choice to use alternative transportation systems is considered. A bilevel mathematical model is proposed considering that the students are reluctant to choose any bus stops that are visited by a route and decide to use an alternative transportation system. In \cite{li2021} a mixed ride approach is proposed, where general and special education students are served on the same bus simultaneously while allowing heterogeneous fleets and mixed loads. Here, mixed rides differ from mixed loads, which refers to buses serving students in different schools using the same bus. In \cite{farzadnia2021} the SBRP is considered a single route service-oriented SBRP. The approach focuses on providing good quality service for students who are sensitive and need more care and safety. In \cite{calvete2021} the SBRP with student choice is introduced. The work is an extension of \cite{calvete2020} to study the student's reaction to the selection of bus stops when they are allowed to choose the bus stop that best suits them. A bilevel optimization model is formulated, considering the selection of the bus stops among the set of potential bus stops and the construction of the bus routes, considering both bus capacity and student preferences.

The ARL strategy attempts to overcome the drawback of the LAR strategy, which tends to create excessive routes. First, students are assigned to clusters while meeting vehicle capacity constraints. Subsequently, the bus stops are selected, and a route is generated for each cluster. Finally, students in a cluster (route) are assigned to a bus stop that satisfies all the problem requirements. The main drawback of this strategy is that the route length cannot be explicitly controlled because each route is generated after the allocation phase and depends on the dispersion of the student within the individual clusters (\cite{Bowerman1995}). The main papers that use a heuristic approach based on the ARL strategy are \cite{chapleau1985, Bowerman1995}. In \cite{chapleau1985}, as in \cite{dulac1980}, multiple schools are considered, but the routes for each school are determined independently. In \cite{Bowerman1995} a multi-objective optimization approach is introduced, which considers efficiency, effectiveness, and equity criteria. The objective is to minimize the number of routes and the total length of the bus route (efficiency), the distance from the student's walk (effectiveness) and optimize both load and length balancing (equity).

Bus stop selection and bus route generation are highly interrelated problems. They are treated separately in previous works due to the complexity and size of the problems, but this generates suboptimal solutions (\cite{salhi1989}). For this reason, recent studies have attempted to solve the problem through an integrated strategy (\cite{Riera-Ledesma2012, Riera-Ledesma2013, Schittekat2013, Kinable2014}). In \cite{Riera-Ledesma2012} the multi-vehicle purchase problem (MV-TPP) is introduced. The MV-TPP is a generalization of the vehicle routing problem and models a family of routing problems that combine stop selection and bus route generation. The goal of the problem is to assign each user to a potential stop to find the least cost routes and choose a vehicle to serve each route so that a route serves each stop with assigned users. Later, an extension of this work is proposed in \cite{Riera-Ledesma2013}. Here, additional resource constraints on each bus route are considered: limits on the distances traveled by students, the number of bus stops visited, and the minimum number of students that a vehicle has to pick up. In \cite{Schittekat2013} a matheuristic to solve large instances is proposed. The matheuristic consists of a construction phase based on a GRASP metaheuristic and an improvement phase based on a VND metaheuristic. The student allocation sub-problem is solved by modeling the problem as a transportation problem. In \cite{Kinable2014} a branch-and-price framework based on a set-covering formulation is proposed. The location of students is not provided, but only the set of stops to which this student can be assigned is known. In \cite{calvete2020} the same problem configuration as in \cite{Schittekat2013} is considered. But there exist two main differences in the solution approach: i) in the construction phase, it is proposed to integrate the allocation of the students and the route construction, while in \cite{Schittekat2013} the integrated sub-problems are the bus stop selection and the route generation, and ii) the procedure includes the solution of four purpose-built Mixed Integer Linear Programming (MILP) problems for selecting the bus stops, together with a more general random selection, which contributes to diversifying the considered solutions. Finally, in \cite{zeng2019covering}, the authors introduced the covering path problem on a grid, which is a bi-objective problem that captures several aspects of SBRP and propose simple techniques to generate feasible solutions.

\subsection{Pricing demand management for routing problems}

Demand management focuses on managing the customer demand to maximize the profitability of a given supply. Demand management can improve profits in two ways \citep{wassmuth2023}: i) increasing revenues by prioritizing high-value customers or serving more customers due to better capacity utilization, and ii) reducing costs by facilitating more eﬃcient order delivery. In addition to proﬁt maximization, demand management can contribute to other goals, such as prioritizing speciﬁc customer groups when demand exceeds capacity or steering customers toward more sustainable delivery times. For a review of demand management, we refer the reader to \cite{fleckenstein2023, wassmuth2023}.


In \cite{campbell2006} the use of incentives to inﬂuence consumer behavior to select time slots to reduce delivery costs is studied. Optimization models for two forms of incentives are proposed and evaluated through simulation studies. In \cite{asdemir2009}, a Markov decision process-based pricing model is proposed. The model recognizes the need to balance the utilization of delivery capacity by the grocer and the need to have goods delivered at the most convenient time for the customer. The model dynamically adjusts delivery prices as customers arrive and make choices. In \cite{yang2016} dynamic pricing policies based on a multinomial logit customer choice model are proposed to determine the incentive (discount or charge) to offer for each time slot when a customer intends to make a booking. In \cite{yang2017}, an approximate dynamic programming method is designed to estimate the opportunity cost of having a customer in a specific area book delivery in a speciﬁc time slot. A demand model, including a multinomial logit model of customers’ delivery time slot choice, is estimated. In \cite{klein2019} differentiated time slot pricing under routing considerations is studied. A mixed-integer linear programming formulation is proposed, in which delivery costs are anticipated by explicitly incorporating routing constraints, and the customer behavior is modeled by a general non-parametric rank-based choice model. In \cite{koch2020} a route-based approximate dynamic programming approach is proposed to influence customer bookings using time window pricing. In \cite{vinsensius2020} a mechanism is proposed to integrate dynamic time slot incentives and order delivery. In \cite{ulmer2020} an anticipatory pricing and routing policy (APRP) method that incentivizes customers to select delivery deadline options efﬁciently for the ﬂeet to fulfill is presented. The dynamic pricing and routing problem is modeled as a Markov decision process (MDP). In \cite{strauss2021} a dynamic pricing approach for standard and ﬂexible time slots is proposed. Flexible time slots are defined as any combination of regular time windows. In \cite{klein2023}, a demand management approach is proposed that explicitly optimizes the combination of delivery spans and prices, which are presented to each incoming customer request. The approach includes an approximation of the value based on an anticipatory sample scenario incorporating a direct online tour planning heuristic.

\section{Mathematical formulations} \label{formulation}

We address the challenge faced by a single school that provides a predetermined incentive for students to opt out of using the school bus. Our focus is on devising a pricing strategy to determine the optimal value of a monetary incentive, denoted as $\tau$. The objective is to ensure that the school experiences an overall reduction in transportation costs. In other words, the savings resulting from operational cost reductions due to a decrease in the number of students using the system must outweigh the expense of providing the incentive to those students choosing to opt out.

Given a set of students using the system, the operational costs for the school are determined by the size of the fleet required to transport all the students and the efficiency of the routes. In addition to transportation constraints, the school must ensure operational conditions to guarantee the safety and comfort of the students. The problem~\ref{routingmodel} presents the general scheme of our cost minimization model.

\begin{align}
\min \quad & \mathcal{L}\left(\text{number of buses ; total travel time}\right)  \label{general.obj1} \\
\text{s.t.} \quad & \text{routing constraints} \label{general.route}\\
&  \text{allocation constraints} \label{general.allocation}\\
& \mathbb{P}\left(\text{overcrowding bus}\right) \leq \alpha && \forall \ \text{bus} \label{general.cap}\\
& \mathbb{E}\left(\text{maximum ride time}\right) \leq t^* && \forall \ \text{bus} \label{general.ride} \\
& \text{walking distance to stop} \leq D && \forall \ \text{students} \label{general.walk}
\end{align}

The primary objective function \eqref{general.obj1} is to minimize the total number of buses required to transport all students. The secondary objective function aims to minimize the total travel time of the route plan, thus encouraging the generation of efficient routes. Constraint \eqref{general.route} ensures the feasibility of routing generation, while constraint \eqref{general.allocation} determines the allocation of students to the selected stops. Constraints \eqref{general.cap} establish an upper limit on the likelihood of overcrowding of buses. The constraint \eqref{general.ride} imposes an upper limit on the maximum expected ride time for a student on any bus. Finally, the constraint \eqref{general.walk} ensures a maximum walking distance from the students' location to their assigned stops.

To formulate the formal mathematical model, we introduce the necessary notation. Let $A = {1, \ldots, N}$ be the set of all students. For each student $i \in A$, the variable $Y_i$ is equal to 1 if the student $i$ opts out of riding the bus and 0 otherwise. We assume that variables $Y_i$ are independent random variables following a Bernoulli distribution with parameter $\theta_i$, representing the opt-out probability for the student $i$. Let $\mathcal{R}\left(\forall i\in A : Y_i=0\right)$ denote the operational cost incurred by the school to transport all students who do not opt out of the system. Note that $\mathcal{R}$ is also a random variable since it is a function of the random variables $Y_i$. 

Then, the school problem is to determine the value of the incentive $\tau$ to minimize the total expected transportation cost, which is given by
\begin{equation}
 \min_{\tau} \quad \mathbb{E}\left[\sum_{i\in A} \tau Y_i + \mathcal{R}\left(\forall i\in A : Y_i=0\right)\right]
\end{equation}



\subsection{Bus stop location problem}

In this section, our focus is on determining the locations of the bus stops and assigning students to these stops.

Students are expected to walk from their homes to their designated bus stops. The school district has policies governing this process, including a maximum allowable walking distance and a maximum number of students per stop. Our main goal is to identify the minimum number of bus stops that comply with all district policies. However, even after determining this minimum, there may be multiple combinations of bus stop locations that result in the same total number of stops. To resolve these ties, we introduce a second objective that minimizes the overall walking distance of the students.

Let $P$ and $A$ represent the set of potential stops and the set of students, respectively. Let $d_{ij}$ be the walking distance from $i\in A$ to $j\in P$, $D$ the maximum distance a student can walk, $p$ the maximum number of students a stop can be assigned. Let $P\left(i\right)=\{j \in P : d_{ij}\leq\delta\}$ be the set of stops within reach of student $i\in A$, and $A\left(j\right)=\{i \in A : d_{ij}\leq\delta\}$ the set of students within reach of stop $j\in P$. Let $z_j$ be a binary decision variable equal to one when a potential stop $j\in P$ is chosen and $y_{ij}$ be a binary decision variable equal to one when the student $i\in A$ is assigned to stop $j\in P\left(i\right)$. 
Then, the allocation model reads as follows:

\refstepcounter{problemno}\label{locationmodel}
\begin{align}
\newProblem :
\min  & \sum_{j\in P} z_j \label{locationmodel:obj1} \\
\min  & \sum_{i\in A}\sum_{j\in P\left(i\right)} d_{ij}y_{ij} \label{locationmodel:obj2} \\
\text{s.t.}  
& \sum_{j\in P\left(i\right)} y_{ij} = 1 && i \in A \label{locationmodel:co1} \\
& y_{ij} \leq  z_j  && i \in A ,\ j\in P\left(i\right) \label{locationmodel:co2} \\
& \sum_{i\in A\left(j\right)} y_{ij} \leq p  && j\in P \label{locationmodel:co3} \\
& y_{ij} \in \{0,1\} && i \in A ,\ j\in P\left(i\right) \label{locationmodel:co4} \\
& z_j \in \{0,1\} && j\in P
\end{align}

where \eqref{locationmodel:obj1} minimizes the number of bus stops for the school and \eqref{locationmodel:obj2} minimizes the total distance students would walk. Constraints \eqref{locationmodel:co1} ensure that every student is assigned one and only one bus stop, constraints \eqref{locationmodel:co2} allow students to be assigned a bus stop that has been chosen as such, and constraints \eqref{locationmodel:co3} limit the number of students per stop.

\subsection{Overbooking} \label{sec:overbooking}


For every student using the bus system, let $X_i$ be a random variable with a Bernoulli distribution $X_i \sim \text{Bernoulli}(\rho_i)$ where $\rho_i$ is the individual ridership of the student (i.e. the probability that the student will use their assigned bus on a certain day). Let $A^*(j)$ be the set of students assigned to stop $j\in P$ (that is, $A^*(j) = \sum_{i\in A(j)} y_{ij}^*$, where $y^*$ is the optimal solution to problem \problemref{locationmodel}). Let $w_{kj}$ be a decision variable that takes value 1 if the stop $j$ is assigned to the bus $k$, and 0 otherwise. Let us further define $Q_k = \sum_{j\in P}\sum_{i\in A^*(j)} w_{kj}X_i$ as the random variable that indicates the true number of students using the bus on a particular day. Then, $Q_k$ is a random variable with a Poisson-binomial distribution with parameters $\rho_i$ for all $i$ such that $i\in A^*(j)$ and $w_{kj} = 1$.


Establishing a linear chance constraint for $Q_k$ is problematic. Therefore, we approximate $Q_k$ with a simpler distribution. In particular, a Poisson-binomial distribution can be approximated for a binomial distribution with the same mean and variance by averaging the parameters $\rho_i$ (see for instance \citep{hoeffding1956distribution, darroch1964distribution}), namely
\begin{equation}
    \mu_{Q_k} = \sum_{j\in P}\sum_{i\in A^*(j)} w_{kj}\rho_i, \qquad \mbox{and} \qquad \sigma_{Q_k}^2 = \sum_{j\in P}\sum_{i\in A^*(j)} w_{kj}\rho_i(1-\rho_i).
    \label{meanvariance}
\end{equation}

Similarly, the binomial distribution can be approximated for normal distributions with the same mean and variance \citep{peizer1968normal}. Therefore, we approximate the random variable $Q_k$ for a normal random variable $\tilde{Q}_k$ with mean $\mu_{Q_k}$ and variance $\sigma_{Q_k}^2$. 

The overbooking constraint for bus $k$ using the approximated normal distribution $\tilde{Q}_k$ reads 
\begin{equation}
    \mathbb{P}(\tilde{Q}_k > Q) \le \alpha, \qquad \forall k\in B,
    \label{chanceconstraint}
\end{equation}
where $Q$ is the maximum capacity of a bus, and $\alpha >0$ is an upper bound for the desired overcrowding probability.

The following proposition establishes linear constraints to reformulate the chance constraint \eqref{chanceconstraint} into a linear programming model. 

\begin{proposition}
For all $k\in B$, the constraints
\begin{align}
& \mu_{Q_k}+\Phi^{-1}\left(1-\alpha\right) \tilde{\sigma}_{\tilde{Q}_k} \leq Q + \frac{1}{2} \label{prop-co01} \\
& \sum_{v=1}^{v^+} v^2 \zeta_v^k \geq \sigma_{Q_k}^2 \label{prop-co02} \\
& \sum_{v=1}^{v^+} v \zeta_v^k = \tilde{\sigma}_{Q_k} \label{prop-co03} \\
& \sum_{v=1}^{v^+} v = 1 \label{prop-co04} 
\end{align}
are valid inequalities for \eqref{chanceconstraint}, where $\zeta_v^k$ is a binary variable and $v^+$ is the maximum possible integer value for $\sigma_{\tilde{Q}_k}$.
\end{proposition}

\subsection{Bus routing problem}

In this subsection, we present the single-school bus routing problem. Let us denote by $D$, $P$, and $S$ the set of depots, stops, and the single school, respectively. We assume that the sets are disjoint and $L:=D\cup P \cup S$ is the set of all locations. Let $\mu_{T_{ij}}$ be the expected value of the travel time between locations $i$ and $j$ where $\left(i,j\right)\in L^2$, $\mu_{T_{i}}$ the expected value of the waiting time or delay at location $i\in P$ and $w_i$ the number of students assigned to stop $i\in P$. Let $B$ be the set of buses and $b_{ik}$ be equal to 1 if bus $k\in B$ starts from depot $i\in D$, and 0 otherwise.

Let $x_{ijk}$ be a binary decision variable equal to 1 when the edge $\left(i,j\right)\in L^2$ is covered by the bus $k\in B$ and 0 otherwise.  Let $w_{kj}$ be a binary decision variable equal to 1 if stop $j\in P$ is assigned to the route $k$, and 0 otherwise. 

The single-school routing problem reads as follows:
\refstepcounter{problemno}\label{routingmodel}
\begin{align}
\newProblem :  \min & \sum_{k\in B} \sum_{i\in D}\sum_{j\in P} x_{ijk} \label{mo4:obj}\\
\min  & \sum_{k\in B} \sum_{i\in L}\sum_{j\in L} \left(\mu_{T_{ij}}+\mu_{T_{i}}\right) x_{ijk} \label{mo4:obj2}\\
\text{s.t. }& \eqref{prop-co01}-\eqref{prop-co04} \\
&   \sum_{i\in D \cup P} x_{ijk} \le w_{kj}, \qquad j \in P \label{mo4:001}\\
& \sum_{k\in B} \sum_{i\in D \cup P} x_{ijk} \le 1 ,\qquad j \in P \label{mo4:002}\\
& \sum_{k\in B} \sum_{j\in P \cup S} x_{ijk}  \le 1 ,\qquad i \in P \label{mo4:003}\\
& \sum_{k\in B} \sum_{i\in L} \left( x_{iik} +  \sum_{j\in D} x_{ijk} +\sum_{j\in S}  x_{jik} \right) =0 \label{mo4:022}\\
& \sum_{i \in D\cup P}  x_{ijk} = \sum_{i \in P \cup S }  x_{jik} ,\qquad  k \in B, j \in P \label{mo4:004}\\
& \sum_{i \in D \cup P}  x_{ijk} \leq \sum_{i \in D}\sum_{j\prime \in P} x_{ij\prime k} ,\qquad k \in B, j\in P \label{mo4:008}\\
& \sum_{j \in L} x_{ijk} \leq b_{ik} ,\qquad k\in B, i \in D \label{mo4:businitialposition}\\
& \mu_{\mathcal{T}_k} - \sum_{i\in D}\sum_{j\in A} \mu_{T_{ij}}x_{ijk}\leq \Delta t_{\text{max}} \qquad \forall k \in B \label{maxridetime} \\
& 1 \leq u_{ik} \leq |P|+2 ,\qquad k \in B , i \in L \label{mo4:st_elim2}\\
& u_{ik} - u_{jk} +\left(|P|+2\right) x_{ijk} \leq |P|+1 ,\qquad  k \in B , i \in L ,j \in L \label{mo4:st_elim4}  \\
\end{align}
Here, the objective function \eqref{mo4:obj} minimizes the number of buses needed, and the objective function \eqref{mo4:obj2} minimize the overall length of the routes. The constraints ensure conditions as follows: Constraints \eqref{prop-co01}-\eqref{prop-co04} ensure that the probability of overcrowding the bus remains less than a certain threshold $\alpha\in (0,1)$. \eqref{mo4:001} implies that a bus travels to the stop only if the stop is assigned to the route of the bus, \eqref{mo4:002} and \eqref{mo4:003} at most one bus arrives and leaves at every stop, \eqref{mo4:022} no bus stays at the same location, nor arrives at a depot or departs from a school, \eqref{mo4:004} is flow conservation for every bus, \eqref{mo4:008} a location can be visited by a bus only if that bus leaves the depot, \eqref{mo4:businitialposition} all buses start their route at their corresponding depot. Constraint \eqref{maxridetime} establishes an upper limit in the expected riding time of students on a bus. Finally, constraints \eqref{mo4:st_elim2} and \eqref{mo4:st_elim4} are the sub-tour elimination constraints where $|P|$ is the number of stops.

Compared to the model in \cite{Caceres2017}, our model incorporates two main novelties. First, we consider a lexicographical objective function that first minimizes the number of buses and then minimizes the total traveling time so the routes become more efficient. Second, we modify the set of constraints \eqref{general.cap} to account for individual ridership, instead of a general ridership per school. Finally, we remark that our model is single-school and therefore several constraints are simplified compared to the multi-school approach of Caceres et al.

\section{Simulation-based solution method} \label{solution}

In this section, we present a simulation-based approach for estimating the expected savings associated with the open-offer policy. It is important to note that the number of students who will choose to reject the bus system when the school offers an incentive $\tau$ is uncertain. To deal with this, we will rely on the concept if {\it opt-out probability}. In particular, the opt-out probability is different for every student and depends on several factors associated with the student's condition, and also on the incentive $\tau$ offered by the school. 

In our approach, we randomly assign decisions to all students based on their respective opt-out probabilities. Subsequently, we establish bus stop locations and routes for the remaining group by excluding the students who have opted out. For this, we solve the allocation problem \problemref{locationmodel} for the remaining students. Then, to address the routing problem described in \problemref{routingmodel}, we employ a decomposition technique that uses column generation, as outlined in \cite[Section 4.3]{Caceres2017}. 

At this stage, we have information on the resulting number of routes and the count of students who have opted out. The potential savings can be determined by calculating the difference between the base routing cost and the cost associated with the particular scenario. By repeating this process, we can generate a wide range of outcomes and approximate the expected savings as the empirical average of the savings of each scenario. 

The following scheme briefly describes our simulation-based approach.

\begin{enumerate}[nosep]
	\item Create sets of students for simulation.
        \item Generate ridership $\bar r$ for the school.
	\item Determine the individual ride for each student.
	\item Determine the individual opt-out probability.
	\item Sample opt-out students and remove them from the system.
	\item Choose stop assignment for the remaining students.
	\item Create a minimal number of routes to solve the problem~\problemref{routingmodel}.
    \item Compute the savings and maximum incentive per student for the scenario.
\end{enumerate}

\subsection{Individual ridership estimation}

In this section, we propose a simple mechanism to estimate the individual ridership level of students. These parameters are then used to create the chance constraint associated with overbooking, as described in Section~\ref{sec:overbooking}. 

In \cite{Caceres2017}, the estimation of student ridership is the average ridership per school, i.e., students from the same school have the same estimate. When considering the fact that students can opt out, the former assumption becomes too simplistic. In particular, if students opt out based on their likelihood of using the bus, then students who continue to use the bus cannot have the same average ridership. To avoid underestimating ridership, we assign each student a ridership estimate that varies depending on how far each one lives from their school.

Let $\rho_i$ be the individual ridership estimate for student $i\in A$ and $\bar{\rho}$ the observed average ridership for one school. 
We aim to find a function $f(d_i)$, where $d_i$ is the distance from the student to the school, which estimates $\rho_i$ so that the average of the resulting estimates matches the original ridership of the school, that is, $\frac{\sum f(d_i)}{n}=\bar{\rho}$. To do so, we assume that $f(d_i)$ is a linear function of $d_i$, namely $\rho_i=f(d_i)=\hat{\rho}_0+\hat{\rho}_1 \times d_i$  where $\hat{\rho}_0$ and $\hat{\rho}_1$ are the tuning parameters to be adjusted. In particular, we assume that $\hat{\rho}_1>0$, so $\rho_i$ increases with $d_i$. Moreover, $d_i$ can be replaced by any increasing function $g(d_i)$.

We solve the following optimization problem to find the coefficients $\hat{\rho}_0$ and $\hat{\rho}_1$ 
\refstepcounter{problemno}\label{ridershipmodel}
\begin{align}
\newProblem : \quad \max_{\{\hat{\rho}_0 , \hat{\rho}_1\}} \quad & \hat{\rho}_1 \label{ridershipmodel:objective} \\
\text{s.t.} \quad & \rho_i=\hat{\rho}_0+\hat{\rho}_1 \cdot g(d_i) && i \in A \label{ridershipmodel:co1} \\
& \sum_{i} \rho_i = n\bar{\rho} \label{ridershipmodel:co2} \\
& 0 \leq \rho_i \leq 1 && i \in A \label{ridershipmodel:co3} 
\end{align}
Here, \eqref{ridershipmodel:objective} aims to maximize the contribution of $d_i$ to the explanation of $\rho_i$. Constraint \eqref{ridershipmodel:co1} defines the linear relation for individual student enrollment. The constraint \eqref{ridershipmodel:co2} forces the average among all individual ridership to match the overall level of overall ridership observed for the school. The restriction \eqref{ridershipmodel:co3} ensures that the individual passengers are valid probability values.

The following proposition provides an explicit formula to compute the individual ridership level.

\begin{proposition}
The problem \problemref{ridershipmodel} admits the solution $\hat{\rho}_1 = \min\left\{\frac{\bar{\rho}}{\bar{d}-\min_i g(d_i) - \bar{d}}, \frac{1-\bar{\rho}}{\max_i g(d_i)-\bar{d}} \right\}$, and $\hat{\rho}_0 = \bar{\rho} + \bar{d}\min\left\{\frac{\bar{\rho}}{\bar{d}-\min_i g(d_i) - \bar{d}}, \frac{1-\bar{\rho}}{\max_i g(d_i)-\bar{d}} \right\}$, where $\bar{d} = \frac{1}{n}\sum_i g(d_i)$. 
\end{proposition}

\begin{proof}
    First, we note that the constraints \eqref{ridershipmodel:co1} and \eqref{ridershipmodel:co2} imply $\hat{\rho}_0 = \bar{\rho} - \hat{\rho}_1\bar d$. Now, as $\hat{\rho}_1$ is non-negative, the constraints \eqref{ridershipmodel:co3} can be replaced by the following two constraints $$0 \le \hat{\rho}_0 + \hat{\rho}_1 \min_i g(d_i), \qquad \hat{\rho}_0 + \hat{\rho}_1 \max_i g(d_i) \le 1.$$
    Replacing the value of $\hat{\rho}_0$ and reordering, we find that 
    $$\hat{\rho}_1(\bar{d}-\min_i g(d_i)) \le \bar{\rho}, \qquad \hat{\rho}_1(\max_i g(d_i) - \bar{d}) \le 1-\bar{\rho}.$$
    Since we are looking for the maximum value for $\hat{\rho}_1$, the result follows. 
\end{proof}

\subsection{Individual opt out probability estimation} \label{sec:individualOOP}

The likelihood of a student opting out of the school bus system can be influenced by various factors, including, but not limited to, the student's distance from the school, family socioeconomic status, parental work schedules, and individual preferences such as safety and environmental considerations, or peer influences. Importantly, the incentive amount offered by the school is a critical determinant in a student's choice to opt out of the bus service.

In our research, we employ a logistic function to determine the likelihood of individual students opting out. The logistic function, commonly utilized for probability modeling \cite{jordan1995logistic}, is defined as
\begin{equation}
    \theta_i(d_i, \tau) = \frac{1}{1+\exp(a\cdot d_i + b\cdot\tau + c)}, \qquad \forall i\in A
    \label{logistic}
\end{equation}
where $\theta_i$ is the opt out probability of student $i\in A$, $d_i$ is the distance to the school, $\tau$ is the incentive offered by the school, and $a, b, c \in \mathbb{R}$ are tuning parameters. 

Parameter tuning for the logistic function \eqref{logistic} can be achieved through a data-driven approach if data is accessible. In the event of data unavailability, one may seek a parameter combination that corresponds to a predetermined expected student behavior based on the school's beliefs. It is commonly understood that the likelihood of opting out decreases as the distance from the school increases. Conversely, the likelihood of opting out is expected to rise with an increase in the incentive parameter $\tau$, indicating that $a > 0$ and $b < 0$. Additionally, as the incentive approaches zero, the likelihood of opting out should tend to zero regardless of the distance from the school. Hence, it is expected that the coefficient $c$ would be positive.

\section{Computational experiments} \label{experiments}

\subsection{Synthetic instances} \label{sec:synthetic}

To create artificial groups of students for the simulation, we take into account three key factors to introduce diversity. Initially, we focus on various types of geographical distribution, including students residing near the school, those living far away, and a uniform distribution of students around the school. To achieve this, we use a beta distribution with parameters $\alpha$ and $\beta$, randomly assigning coordinates within a square area with a side length of 3 miles. The coordinates of each student are generated independently with the following parameters: $\left(\alpha=1,\ \beta=1\right)$, $\left(\alpha=0.45,\ \beta=0.45\right)$, and $\left(\alpha=6.5,\ \beta=6.5\right)$. Second, we acknowledge the significance of considering the school's student population as another variable. To thoroughly investigate the influence of student numbers on the simulation results, we carried out experiments with two distinct sample sizes: one consisting of 400 students and another with 800 students. Lastly, we examine the average ridership per school. To gain a comprehensive understanding of student transportation needs, we specifically analyze scenarios with an average ridership of $\bar{r} = 0.3$ and $\bar{r} = 0.8$. Through intentional manipulation of three key factors, geographical spread, student enrollment, and average ridership, we successfully generated a variety of artificial students for simulation. 



\begin{figure}[ht]
	\begin{center}
		\includegraphics[width=1\columnwidth]{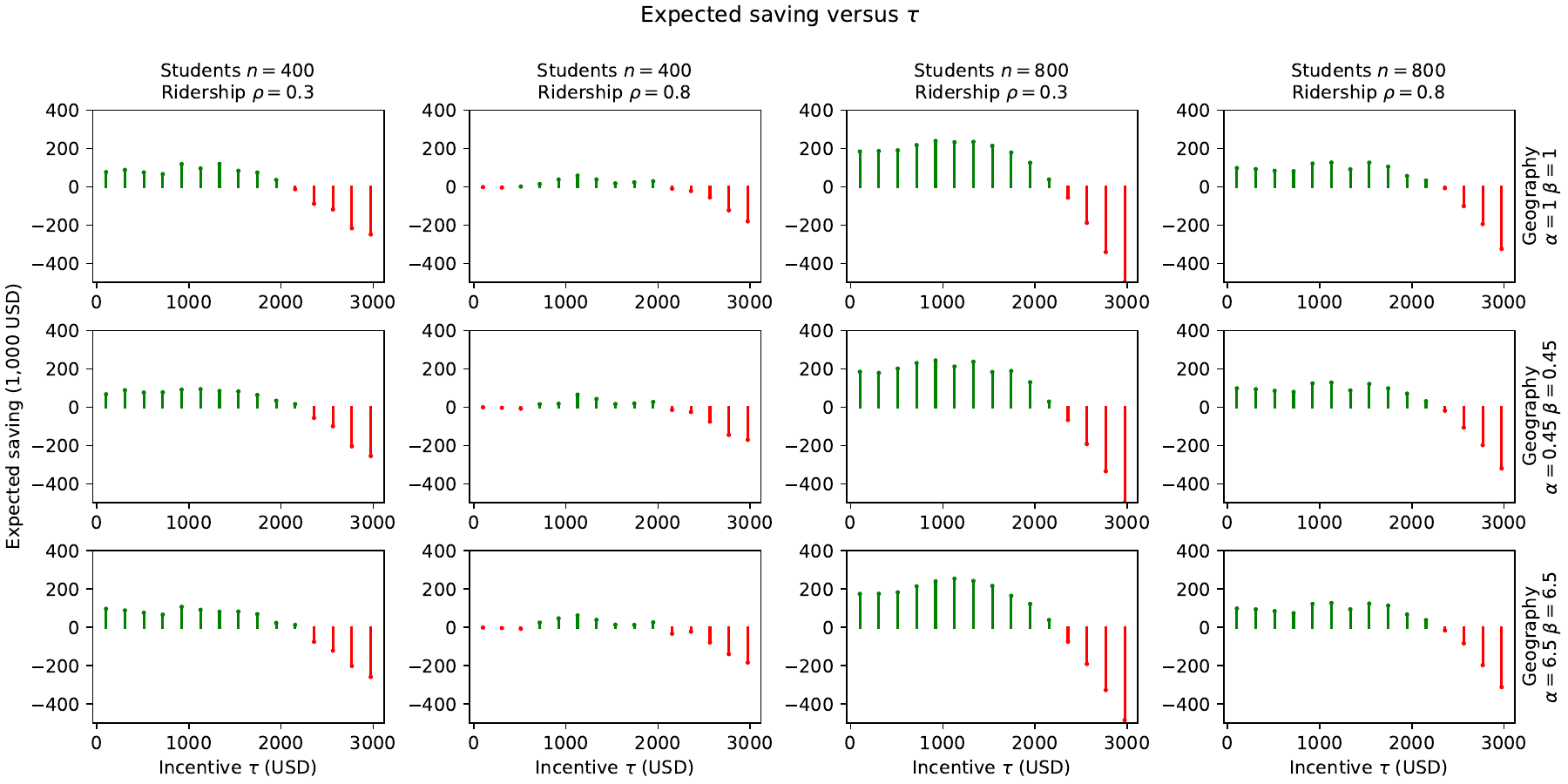}
		\caption{Expected savings for each data set} \label{fig:pic_expected_saving_set_all}
	\end{center}
\end{figure}

Figure~\ref{fig:pic_expected_saving_set_all} shows the expected savings in 12 data sets. We performed 100 repetitions of the simulation procedure for each data set and incentive value ($\tau$), then averaged the results to estimate the expected savings. Based on our analysis, we find that the open-offer policy offers significant benefits, particularly for schools with a high student population and low ridership. Even modest incentives can lead to a notable reduction in the bus fleet, with potential savings of up to two buses for schools with 400 students and up to three buses for schools with 800 students.

However, we also observe that as the incentive size increases, the overall effectiveness of the policy decreases. This is mainly due to the fact that while the incentive amount increases, there is no corresponding increase in the number of students opting out of the bus system. This phenomenon indicates that there is a threshold beyond which offering larger incentives may not produce substantial reductions in fleet size. Therefore, careful consideration of the incentive amount is crucial to ensure that the policy remains cost-effective and efficient for the school district.

\section{Case study: the Williamsville Central School District } \label{sec:casestudy}

The Williamsville Central School District (WCSD) is an esteemed public school district located in Williamsville, New York, United States. The WCSD serves the suburban communities of Williamsville, Amherst, and Clarence in Erie County, comprising a total of 13 schools in the district. It includes six elementary schools (EL), four middle schools (MI), and three high schools (HI).  ELs are strategically located throughout the district and provide education from kindergarten through fifth grade. The MIs are for students in grades six through eight, while the HIs are for students in grades nine through twelve.

Encompassing an area of approximately 40 square miles, WCSD accommodates over 10,000 students, nurturing their academic growth and development. The district strongly emphasizes student transportation, recognizing the importance of safe and efficient travel. To facilitate this, WCSD operates a fleet of nearly 100 buses, maintaining a ratio of approximately 2:3 between its own buses and those provided by a contractor. This robust transportation system ensures that students can access the various schools within the district effectively.

\begin{figure}[ht]
    \centering
    \includegraphics[width=0.25\columnwidth]{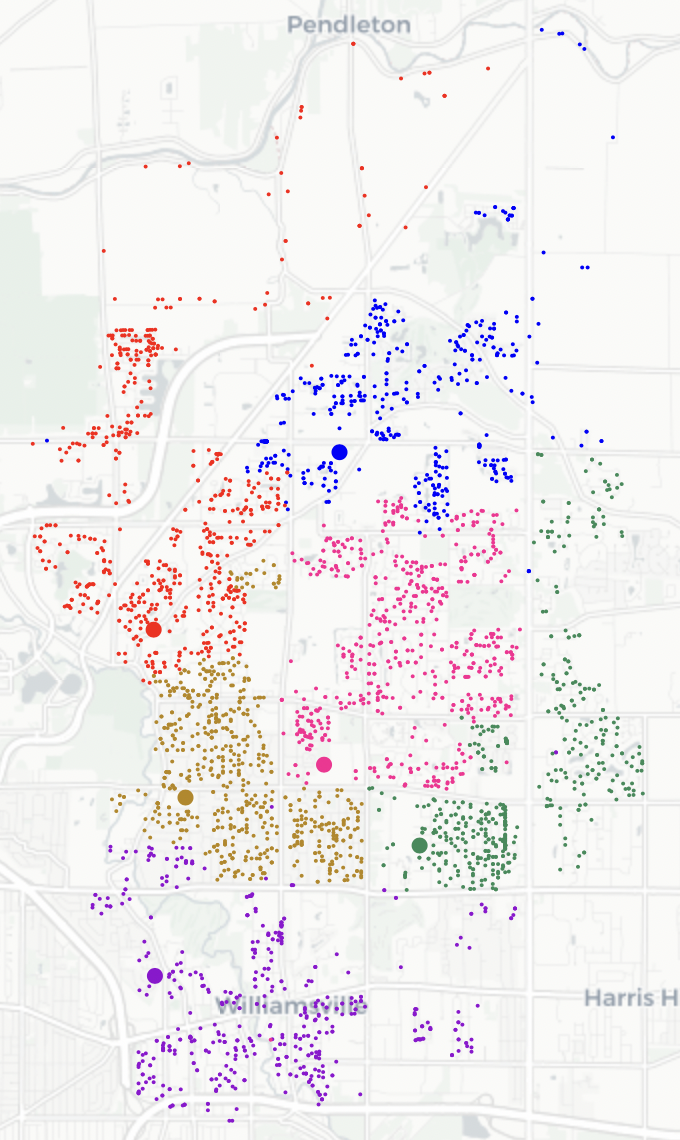}
    \includegraphics[width=0.25\columnwidth]{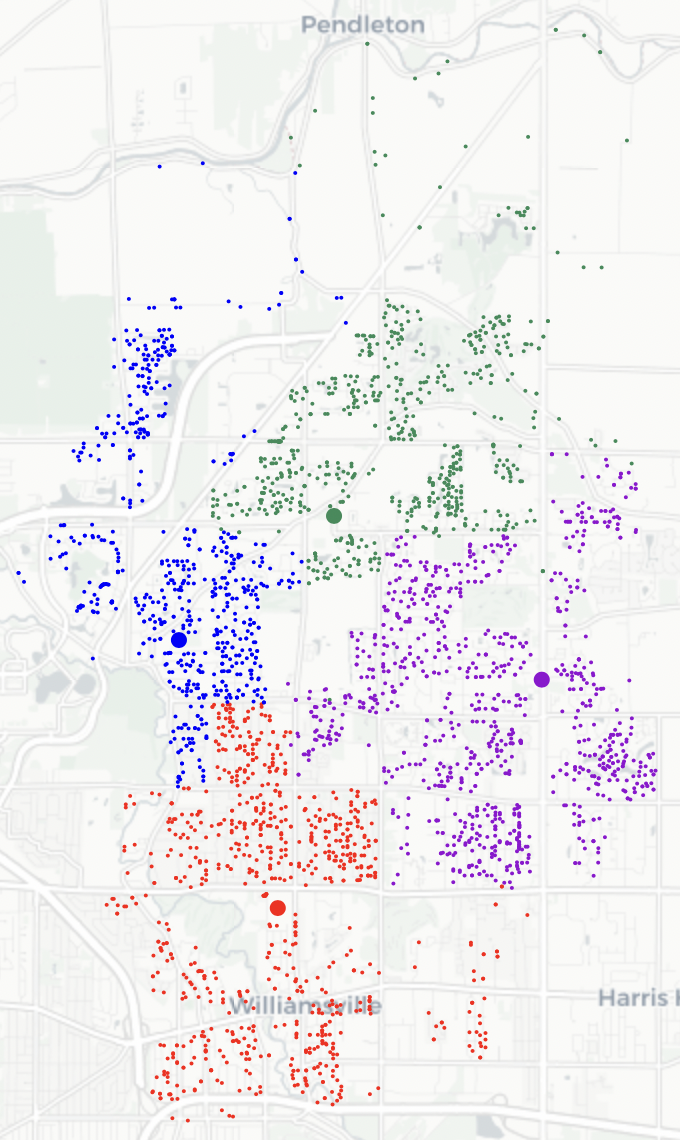}
    \includegraphics[width=0.25\columnwidth]{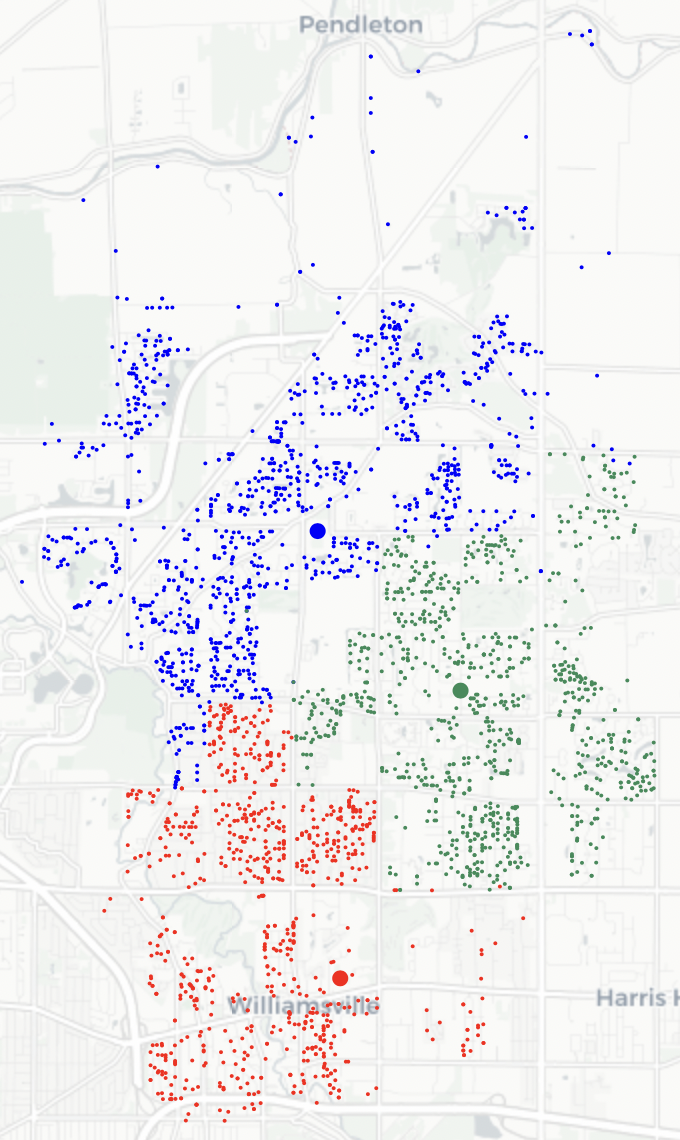}
    \caption{Dispersion for students for the 13 schools in the WCSD in 2014 - 2015: six elementary schools (left); four middle schools (center); three high schools (right).}
    \label{fig:dispersion}
\end{figure}

The spatial distribution of the students is shown in Figure~\ref{fig:dispersion}. Elementary, middle and high schools exhibit comparable spatial patterns. However, the disparity in school numbers results in high school students being widely scattered than those in elementary and middle schools. On the other hand, Figure~\ref{tax:ridership} reveals that high schools have the lowest level of ridership, meaning that in practice, a considerable number of students do not use the bus service.

\begin{figure}[ht]
    \centering
    \includegraphics[width=1\columnwidth]{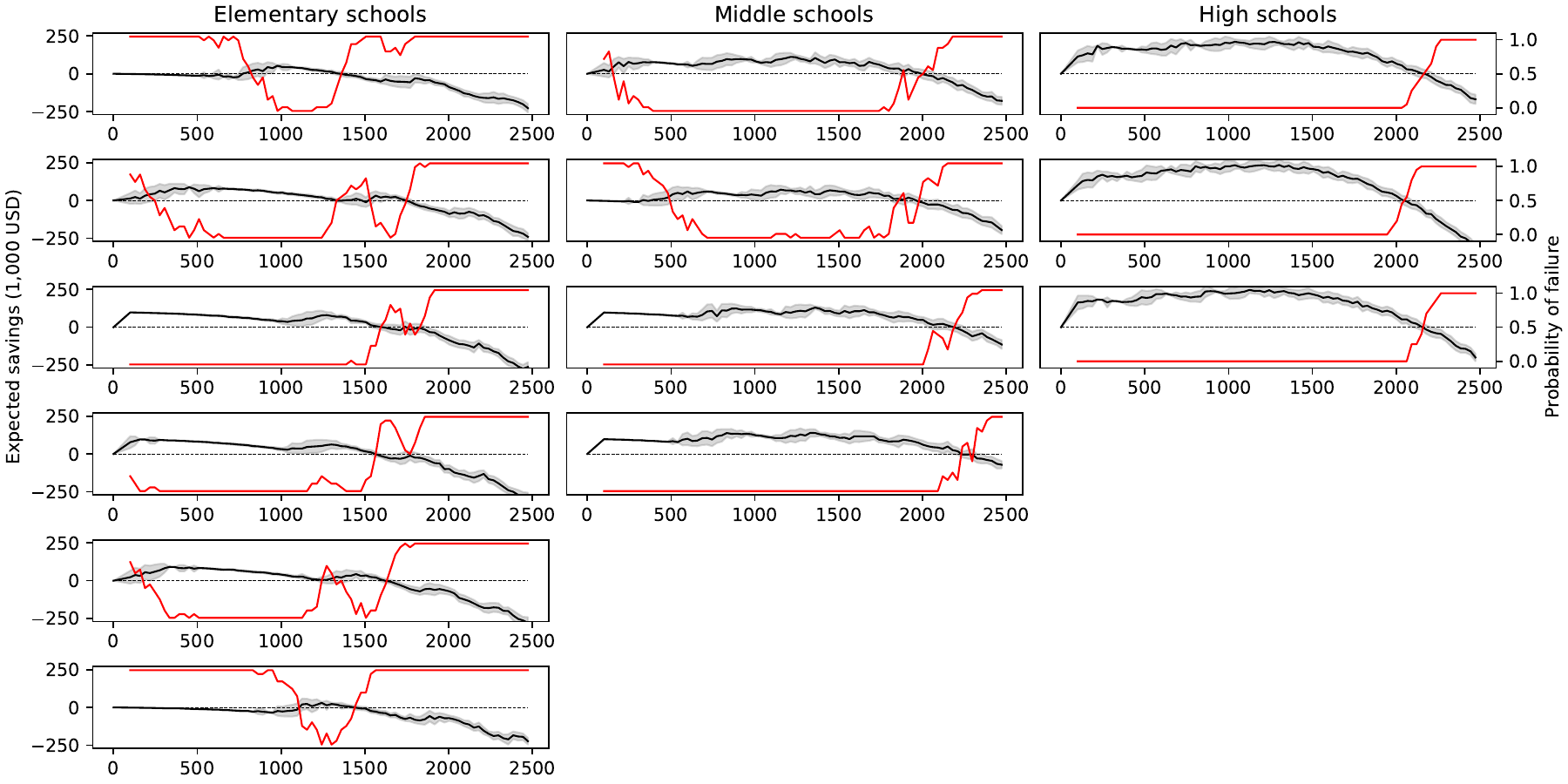}
    \caption{Expected savings (black line) and probability of failure (red line) for each school of the WCSD. The gray area mark the confident interval given by the standard deviation of the expected saving.}
    \label{fig:results}
\end{figure}

Figure~\ref{fig:results} presents a comprehensive overview of the expected savings for each school, highlighting the influence of various incentive values ranging from 0 to 2,500 USD. It also tracks the evolution of the {\it probability of failure}, a metric calculated as the proportion of replicas that results in negative savings.

The data reveal a general trend of unfavorable results for elementary schools under the open-offer policy. In particular, schools 1 and 6 are particularly affected, with consistently negative expected savings in most incentive values. Even when the savings are positive, they remain minimal. Schools 2, 3, 4, and 5 exhibit similar patterns, with rapidly decreasing savings and increasing probabilities of failure as incentives increase. The middle schools exhibit a better behavior, with potential savings of approximately 160,000 USD achieved at incentive levels of around 1,000 and 1,500 USD. Interestingly, schools 7 and 8 face a higher risk when offering low incentives (below 500 USD), unlike schools 9 and 10. This disparity can be attributed to the lower enrollment in schools 8 and 9, which makes students in schools 9 and 10 more inclined to accept lower incentives to opt out of the system. For high schools, the results are more consistent. All three schools display similar behavior, with substantial savings of up to around 260,000 USD observed for incentives up to 1,500 USD. Notably, the probability of failure remains low, only increasing significantly with very high incentive values. Even with lower incentives, the savings are considerable, highlighting the impact of low ridership levels in these schools.

In summary, our findings indicate that the open-offer policy demonstrates promise in middle and high schools, resulting in an average reduction of two and three buses, respectively, for reasonable incentives. However, its application in elementary schools appears less effective, failing to yield a significant reduction in bus usage. Several factors may explain this observation. First, the level of ridership directly influences the number of students who opt out, with higher ridership typically associated with higher opt-out rates. Additionally, bus capacity differs between school levels, with middle and high schools having lower capacities compared to elementary schools. Consequently, a greater proportion of students in elementary schools must opt out to achieve a significant reduction in bus usage.

\section{Conclusions} \label{discussion}

In this paper, we introduce a new variant of the School Bus Routing Problem called School Bus Routing Problem with Open Offer Policy (SBRP-OOP), which focuses on finding a pricing strategy that balances incentive payments for students who choose not to use the bus with the savings achieved by operating fewer buses. After excluding students who accept the incentive from the routing process, our model addresses multiple objectives simultaneously for a single school. This involves selecting the optimal set of stops from a pool of potential stops, assigning students to stops while considering factors like maximum walking distance and student capacity, and generating routes that minimize bus travel distance.

To evaluate the effectiveness of the open offer policy, we performed simulations using both synthetic and real data. Our analysis indicates that this policy is particularly advantageous for schools with a sizable student population and low ridership, where even modest incentives can result in significant savings. However, it is crucial to consider that as incentive sizes increase, the overall effectiveness of the policy diminishes. For real data from Williamsville Central School District (WCSD), we found that the open offer policy can be implemented with risk and promising benefits in certain middle and high schools. However, for some elementary schools, implementing the policy seems impractical.

While these observations suggest possible advantages, it is crucial to exercise care in interpreting these findings due to the influence of various factors like local infrastructure, geographical limitations, public policies, and specific school transport regulations on the results. Consequently, further investigation and examination are necessary to confirm and investigate these results in a range of different scenarios.

\section*{Acknowledgment}
This research was supported in part by the Chilean Science and Technology National Council (CONICYT), grant FONDECYT 11181056. Computational resources were sponsored by the Department of Industrial and Systems Engineering and provided by the Center for Computational Research (CCR) of the University at Buffalo, New York.


\end{document}